\documentclass[a4paper,leqno,twoside]{amsart}
\usepackage[british]{babel}
\usepackage{amscd}
\usepackage{url}
\usepackage{tabularx}

\theoremstyle{plain}
\newtheorem{thm}{Theorem}

\newtheorem{prop}[thm]{Proposition}

\theoremstyle{definition}
\newtheorem{df}{Definition}
\newtheorem{dfs}[df]{Definitions}
\newtheorem{ex}{Example}

\theoremstyle{remark}
\newtheorem*{rmk}{Remark}
\newtheorem*{acks}{Acknowledgments}

\newcommand{\ie}{\textit{i.e. }}
\newcommand{\eg}{\textit{e.g. }}
\newcommand{\p}{\mathbb{P}}

\newcommand{\Gr}{\mathbb{G}}

\newcommand{\OO}{\mathcal{O}}

\newcommand{\gen}[1]{\langle#1\rangle}

\DeclareMathOperator{\Sing}{Sing}

\DeclareMathOperator{\Pfaff}{Pfaff}
\def\cocoa{{\hbox{\rm C\kern-.13em o\kern-.07em C\kern-.13em o\kern-.15em A}}}
\DeclareMathOperator{\PGL}{\mathbb{P}GL}
\DeclareMathOperator{\ii}{II}

\begin{document}

\title{Linear Congruences and hyperbolic Systems of conservation Laws}

\author{Pietro De Poi \and Emilia Mezzetti}

\thanks{This research was partially supported by the DFG Forschungsschwerpunkt
``Globalen Methoden in der Komplexen Geometrie'' for the first
author, by   MIUR, project \lq\lq Geometria sulle variet\`a
algebriche'' for the second author
and by funds of the University of Trieste (fondi
60\%) and the EU (under the EAGER network) for both.}

\address{Mathematisches Institut\\
Universit\"at Bayreuth\\
Lehrstuhl VIII \\
Universit\"atsstra\ss e 30\\
D-95447 Bayreuth\\
Germany\\}
\email{Pietro.DePoi@uni-bayreuth.de}

\address{Dipartimento di Matematica e Informatica\\
Universit\`a degli Stud\^\i\ di Trieste\\
Via Valerio, 12/b\\
I-34127 Trieste\\
Italy}
\email{mezzette@univ.trieste.it}
\keywords{Congruences
of lines, Grassmannian of lines, Conservation Laws, T-systems,
fundamental points}
\subjclass[2000]{Primary 14M15, 35L65  Secondary 53A25, 53B50}

\date{\today}

\begin{abstract}
S. I. Agafonov and E. V. Ferapontov have introduced a construction
that allows naturally associating
to a system of partial differential equations of
conservation laws a congruence of lines in an appropriate
projective space. In particular hyperbolic systems of Temple class
correspond to congruences of lines that place in
planar pencils of lines. The language of Algebraic Geometry turns out
to be very natural in the study of these
systems. In this article, after recalling the
definition and the basic facts on congruences of lines,
Agafonov-Ferapontov's construction is illustrated
and  some results of classification for Temple systems are presented.
In particular, we obtain the classification of
linear congruences in $\p^5$, which correspond to some classes of
$T$-systems in $4$ variables.
\end{abstract}

\maketitle

\section*{Introduction}
Linear congruences of lines in $\p^n$ are the (irreducible)
subvarieties of dimension $n-1$ of the Grassmannian $\Gr(1,n)$,
embedded in $\p^{\binom{n+1}{2}-1}$
by the Pl\"ucker embedding,  obtained by the
intersection  with a linear space of dimension  $\binom{n}{2}$.

Such a congruence of lines ${\mathcal B}$ has, in general, order one, \ie
through a general
point in $\p^n$ there passes
only one line of ${\mathcal B}$. Moreover the family of lines
parametrised by ${\mathcal B}$
can be characterized as the set of the $(n-1)$-secant
lines
of the fundamental locus $\Phi\subset\p^n$, where $\Phi$ is defined
by the property that through a point in
it there pass infinitely many lines of ${\mathcal B}$.

This article deals with a recent discovered application of
congruences of lines to mathematical physics
(precisely to hyperbolic systems of conservation laws)
due to S. I. Agafonov and E. V. Ferapontov: see \cite{AF} and \cite{AF1}.
More precisely, to a system of conservation laws, which has the form
$\frac{\partial u^i}{\partial t}+\frac{\partial f^i({\mathbf
u})}{\partial x}=0$, with
$i=1,\dotsc,n-1$, they associate  an $(n-1)$-parameter family
${\mathcal B}$ of lines
in $\mathbb{P}^n$, defined by the parametric equations
$y_i=u^i\lambda-f^i({\mathbf u})\mu$, $i=1,\dotsc,n-1$ and
$y_0=\lambda$, $y_n=\mu$, where
$(\lambda:\mu)\in\mathbb{P}^1$ are the parameters of a line of
${\mathcal B}$,
${\mathbf u}=(u^1,\dotsc,u^{n-1})$
are the (local) parameters of ${\mathcal B}$ and $(y_0:\dotsb:y_n)$
are the
homogeneous coordinates on
$\mathbb{P}^n$. It turns out that with this correspondence
the basic concepts of the theory of the
systems of conservation laws acquire a clear and simple projective
interpretation.
For instance, for a particular class of systems of conservation laws,
the so called
$T$-systems (see \cite{AF1} and \cite{AF2}), their corresponding
family of lines ${\mathcal B}$ is characterized
by the fact that the lines of ${\mathcal B}$ passing through a  point of its
focal locus form a planar pencil
of lines; moreover, a reciprocal transformation (see \cite{AF2} for a
definition)
of one of these systems of conservation laws corresponds to a
projectivity in $\p^n$, and
vice versa. Therefore, the classification of the $T$-systems
is equivalent to the study of these families of lines ${\mathcal B}$.

Classically, the study of congruences of lines in $\p^3$ was
started by E.~Kummer in \cite{K}, in which he gave
a classification of those of
order one. More recently Z. Ran in \cite{R} studied
the surfaces of order one in a
general Grassmannian $\Gr(r,n)$ \ie families of $r$-planes in $\p^n$
for which the general $(n-r-2)$-plane meets only one element of the family.
He gave a classification of such surfaces, obtaining in particular,
in the case of $\p^3$, a modern and more correct proof of Kummer's
classification.

The congruences of lines in $\p^4$, with special
regard to those of order one, have been considered
by G. Marletta in \cite{M3} and \cite{M1}. The classification of the linear
ones
has been given by G. Castelnuovo in
\cite{C}, where a detailed description of these particular subvarieties of
$\Gr(1,4)$ is obtained. These classical results in $\p^4$ have been
analysed and extended by P. De Poi in \cite{DP2}, \cite{i},
\cite{ii}, and
\cite{iii}.

A general fact about linear congruences in $\p^n$ is that
the lines of the family passing through a general
focus form a linear pencil and the plane
of this pencil cuts the focal locus residually along a plane curve of degree
$n-2$. So linear congruences always define Temple systems.
Conversely for $n\leq 4$ in \cite{AF2} it has been proved that  all families of
$T$-systems are even algebraic, and more precisely linear congruences.

In higher dimensions nothing is  known, in particular also a
complete classification of the linear congruences of lines is still missing.
So in this paper we have initiated a systematic study of the more
simple unknown case, that
of linear congruences in $\p^5$.

For general linear congruences in $\p^5$, the focal locus is a smooth
Palatini threefold, which is a scroll over a cubic surface $S$ in
$\p^3$ (see \cite{O},
\cite{FM}). This surface $S$ can be realized as follows:
let the congruence ${\mathcal B}$ be defined
as $\Gr(1,5)\cap\Delta$, for a $10$-dimensional
linear space $\Delta$. The dual of the Grassmannian is a cubic
hypersurface (the ``Pfaffian'', see
Section~\ref{sec:4}), then $S$ is
naturally identified
with the intersection
of $\check\Gr(1,5)$ with the dual of $\Delta$. Classifying
linear congruences in $\p^5$ amounts to describe
all special positions of the $3$-space $\check\Delta$ with respect to
$\check\Gr(1,5)$ and to its singular locus.
For instance, when $S$ meets $\Sing \check\Gr(1,5)$, the focal locus
acquires some linear
irreducible components.
Particularly interesting are the cases  when $S$  splits: the  description of
these congruences relies on a recent classification of the linear
systems of $6\times 6$
skew-symmetric matrices
of constant rank $4$ up to the natural action of the projective linear group
$\PGL_6$ (\cite{MM}).

This article is structured as follows: in Section~\ref{sec:1}, the
basic definitions connected to congruences of lines in the algebraic
setting are given.
In Section~\ref{sec:2} we recall some definitions and results about
systems of conservation laws, reciprocal transformations and
systems of Temple class.
In Section~\ref{sec:3} the correspondence between systems of
conservation laws and families of lines is illustrated, with
special regard
to systems of conservation laws of Temple class. In
Section~\ref{sec:4} we collect some general facts about linear congruences.
Finally, in Section~\ref{sec:5}
we study the  linear congruences in $\p^5$. We
first consider the case in which the cubic surface $S$ has one or
more singular points on the singular locus of $\check\Gr(1,5)$,
these points correspond to some 3-dimensional linear spaces that
enter in the focal locus. We then study the congruences such that
the surface $S$ is reducible,  getting four types of congruences. In
some cases the focal locus has a parasitic component, \ie an irreducible,
maybe embedded component, which is not met by a general line of $\mathcal
B$.

In a forthcoming paper we plan to apply these results to the
classification of Temple systems
in 4 variables. We have  to point out that the classification
considered here holds
over an algebraically
closed field, so it will be necessary to refine it over the real field.

\begin{acks}
We wish to thank Jenya Ferapontov for
introducing us to this beautiful connection between algebraic
geometry and partial differential equations. We also thank
Dario~Portelli and Giorgio~Tondo for interesting discussions.
\end{acks}

\section{Notation, Definitions and Preliminary Results}\label{sec:1}

In the realm of algebraic geometry,
we will work with schemes and varieties over $\mathbb C$, with
standard notation and definitions as in \cite{H}. A
\emph{variety} will always be projective.
We refer to \cite{DP2} and \cite{D2} for general results and references
about families of lines, focal diagrams and congruences, and to \cite{GH} for
notations about Schubert cycles. In particular we denote  by $\sigma_{a_0,a_1}$
the Schubert cycle of the lines in $\p^n$ contained in a fixed
$(n-a_1)$-dimensional subspace $H\subset\p^n$ and which meet a fixed
$(n-1-a_0)$-dimensional subspace $\Pi\subset H$.

\begin{dfs} A  \emph{congruence of lines ${\mathcal B}$ in $\p^n$}
is a flat family of lines of dimension $n-1$, and we can think of it as
a $(n-1)$-dimensional subvariety of the Grassmannian
$\mathbb G(1,n)$. Its \emph{order} $a_0$
is the number of lines passing through a general point in $\p^n$.
\end{dfs}

Throughout this article, we will denote by $\Lambda:=\{(b,P)\mid P\in
\Lambda(b)\}$
the incidence correspondence associated to the family ${\mathcal B}$,
with its natural projections $f$ and $p$ to $\p^n$ and to ${\mathcal
B}$ respectively.
If $b\in {\mathcal B}$, then $\Lambda(b):=f(\Lambda_b)$ denotes the line
parametrised by $b$ and $\Lambda_{b}:=p^{-1}(b)$
the corresponding subset of $\Lambda$. We can summarize all in the diagram:
\begin{equation}\label{cd:f}
\begin{CD}
\ \ \ {\Lambda_{b}}\subset
@. {\Lambda} @. \subset {\mathcal B}\times \p^n @>{f}>>
\p^n\supset {\Lambda(b)} \\
@VVV @VVpV\\
\ \ \ b\in @. {\mathcal B}.
\end{CD}
\end{equation}

Note that $f$ is surjective if and only if $a_0>0$. In this case $f$
is a  map of degree
$a_0$.

The basic objects when one deals with congruences are the focal and
fundamental loci.
\begin{dfs} Let ${\mathcal B}$ be a congruence of lines; then
the  \emph{focal divisor}
of the family ${\mathcal B}$  is
the ramification divisor ${R}\subset \Lambda$ of
\begin{equation}\label{eq:f}
f:\Lambda\rightarrow \p^n;
\end{equation}
the schematic image $F$ of $R$ under $f$ is called the
\emph{focal locus}:
${F}=f(R)\subset \p^n$. The \emph{foci} of ${\mathcal B}$ are the
branch points of $f$.
The  \emph{fundamental locus}
${\Phi}$ is the set of points $y$ contained in more lines of the
family than expected:
$$
\dim f^{-1}(y)> n-\dim f(\Lambda).
$$
\end{dfs}

\begin{ex}
The family of the tangent lines to a curve $C$ (which is not a line) in $\p^2$
is the most simple example of a congruence.
Other examples are the family of the secant lines
to a curve in $\p^3$ or the one of the tangent lines to
a surface in $\p^4$.
\end{ex}

\begin{rmk}
The fundamental locus $\Phi$ is in general properly contained in the
focal locus $F$,
but if ${\mathcal B}$ is a congruence of order one, then $\Phi=F$ and
the codimension
of $F$ in $\p^n$ is $\geq 2$.
\end{rmk}

The following recent result  gives to converse  statement:

\begin{thm}(F. Catanese, P. De Poi, 2004, \cite{D}) \label{DC}
Let ${\mathcal B}$ be a congruence such that the fundamental locus $\Phi$
coincides (set-theoretically) with the focal locus $F$;
then the order of ${\mathcal B}$ is zero or one.
\end{thm}

\begin{rmk} It is important to note that the focal locus $F$
often has some unexpected components. For example, let ${\mathcal B}$ be the
congruence of the secant lines to a curve $C$ in $\p^3$. If $L$
is a line meeting $C$ at two points $x,y$ such that
the tangent lines to $C$ at $x$ and $y$ are incident, then
the tangent plane to ${\mathcal B}$ at the point corresponding to $L$
     is contained in the Grassmannian
and $L\subset F$. Since a curve in $\p^3$ has in general a one-dimensional
family of secant lines of this type (called stationary secants), the
congruence ${\mathcal B}$
will have a focal surface. The only curve in $\p^3$ without
stationary secant lines
is the twisted cubic. For more details, see \cite{ABT}.

The focal locus of a congruence  may also have a component which
is not met by a general line of $\mathcal B$. Such a component is
called a parasitic component. Some explicit examples will be
described in Sections~\ref{sec:4} and \ref{sec:5}.
\end{rmk}

     From now on, we assume that $f$ is surjective. In this case, the
fundamental locus is the set of points $P\in\p^n$
for which the fibre of the map $f:\Lambda\rightarrow\p^n$ has
positive dimension.

\begin{thm}\label{prop:fofi} (C. Segre, 1888, \cite{Sg}, C.
Ciliberto and E. Sernesi, 1992,
\cite{SC})
With notations as above, given a congruence ${\mathcal B}$, for the general
$b\in {\mathcal B}$, the corresponding line
$\Lambda(b)\subset\p^n$ contains exactly $n-1$ foci (counting
multiplicities) which
are foci for the line $\Lambda(b)$. Otherwise $\Lambda(b)\subset F$.

Moreover, if $\dim(F)=n-1$,  then $\Lambda(b)$ is tangent to $F$ at
its (smooth)
focal non fundamental points.
\end{thm}

Let us now give some interesting and useful examples of first order
congruences.

\begin{itemize}
\item \emph{Congruences of order one in $\p^3$}  (E. Kummer, 1866, \cite{K},
Z. Ran, 1986, \cite{R}): these are all classified and can be divided
in three cases:
\begin{enumerate}
\item $F$ is a skew cubic, ${\mathcal B}$ is the family of its secant lines. In
the Pl\"ucker
embedding, ${\mathcal B}$ is a Veronese surface;
\item $F=C\cup L$, $C$ is a rational curve of degree $d$, $L$ is a
$(d-1)$-secant line of $C$,
${\mathcal B}$ is the family of lines meeting both $C$ and $L$;
\item $F$ is a double structure on a line $L$, ${\mathcal B}$ is a union of
pencils of lines with centres 
at the points of $L$.
\end{enumerate}

\item \emph{Linear congruences in $\p^4$} (G. Castelnuovo, 1891, \cite{C}):
if ${\mathcal B}=\Gr(1,4)\cap \p^6$, then $F$ is a projected Veronese
surface or
a degeneration of
it, ${\mathcal B}$ is the family of the trisecant lines of $F$. The
order of
${\mathcal B}$ is $1$. The lines in ${\mathcal B}$
through a general point of $F$ form a planar pencil and cut a conic on $F$.

\item\emph{$(n-1)$-secant lines of smooth codimension two
subvarieties in $\p^n$} (P. De Poi, 2003, \cite{D2}):
the congruences of order one formed
by these lines are completely described,
they are all contained in $\p^n$, with $n\leq 5$.
\end{itemize}

\section{Systems of Conservation Laws}\label{sec:2}

In the realm of mathematical physics,
we will work over $\mathbb R$. All the functions are---at least---$C^1$.

\begin{dfs} A \emph{system of conservation laws} is a quasi-linear
system of first order partial differential equations of the form
\begin{align}\label{eq:cons}
{\frac{\partial u^i}{\partial t}+\frac{\partial f^i({\mathbf
           u})}{\partial x}}&=0 & i&=1,\dotsc,n-1
\end{align}
where ${\mathbf u}(x,t)=(u^1(x,t),\ldots , u^{n-1}(x,t))$ are the
unknown functions and
the
$f^i({\mathbf u})$'s are functions defined over a domain
${\Omega}\subset{\mathbb R}^{n-1}$.

The system can be written:
\begin{align}\label{eq:cons2}
{u^i_t + \Sigma_i\frac{\partial f^i({\mathbf u})}{\partial u^j}}u^j_x&=
u^i_t+Jf({\mathbf u})\cdot u_x=0& i&=1,\dotsc,n-1
\end{align}
where
$Jf$ denotes the Jacobian matrix of $f$.

The system is called \emph{hyperbolic} (resp.
\emph{strictly hyperbolic})
if all the eigenvalues of $Jf$
are real (resp. real and distinct).
\end{dfs}

\begin{dfs}
If the system \eqref{eq:cons} is strictly hyperbolic, the eigenvalues
$$
{\lambda_1(\mathbf u)}<{\lambda_2(\mathbf u)}<\dotsb<
{\lambda_{n-1}(\mathbf u)}
$$
are called \emph{characteristic velocities}.

The integral trajectories  $\gamma_i$ of  the fields of eigenvectors
${\mathbf v}_i $ are called the
      \emph{rarefaction curves}:
$$
\dot{\gamma}_i(t)={\mathbf v}_i(\gamma_i(t)).
$$
\end{dfs}

\begin{rmk}
Given a strictly hyperbolic system of conservation laws as above,
through a point ${\mathbf u}\in \Omega\subset {\mathbb R}^{n-1}$
there pass $n-1$ rarefaction
curves. Indeed  each eigenvector of the eigenvalue
$\lambda_i(\mathbf u)$, for $i=1,\dotsb, n-1,$ has the direction of the tangent line
to the curve $\gamma_i$.
\end{rmk}


\subsection{Temple Systems} These systems, which are strictly
hyperbolic, were introduced by B. Temple
in \cite{T}. They  naturally arise in the theory of equations of
associativity of
2D topological field theory (see \cite{Du}).

\begin{df} A  strictly hyperbolic system \eqref{eq:cons} is said
to be  \emph{linearly degenerate} if
$$
{L_i(\lambda_i)(\mathbf u)}=0\quad\forall i
$$
where $L_i$ is the {Lie derivative}
in the direction of ${\mathbf v}_i$.
\end{df}

This means that the eigenvalue $\lambda_i$ is constant along the
rarefaction curve
$\gamma_i$.

\begin{df}
A strictly hyperbolic system \eqref{eq:cons} is said
to be a \emph{Temple system} or a \emph{$T$-system} if
it is linearly degenerate and the rarefaction curves are lines in the
coordinates $(u^1,\dotsc,u^{n-1})$ (see \cite{T}).
\end{df}


\subsection{Reciprocal Transformations}

Let us define two new independent variables $(X,T)$ by
\begin{align}\label{eq:rec}
dX&:=\sum_i \alpha_i (u^idx+f^i(\mathbf u)dt)+\mu dx+\nu dt\\
dT&:=\sum_i \tilde{\alpha}_i (u^idx+f^i(\mathbf u)dt)+\
\tilde{\mu}
      dx+\tilde{\nu} dt.
\end{align}
By \eqref{eq:cons} the two $1$-forms on the right are closed, so
\eqref{eq:cons} takes the new form
\begin{align}
{\frac{\partial U^i}{\partial T}+\frac{\partial F^i({\mathbf
           U})}{\partial X}}&=0 & i&=1,\dotsc,n-1
\end{align}

\begin{df}
Transformations of the form \eqref{eq:rec} are called
\emph{reciprocal}.
\end{df}

Reciprocal transformations are known to preserve the class of $T$-systems
(see \cite{AF}).


\section{The Correspondence}\label{sec:3}

Agafonov and Ferapontov ((\cite{AF} and \cite{AF1})) associate to a
system of conservation
laws \eqref{eq:cons}  a {$(n-1)$-parameter family ${\mathcal B}$ of
lines in $\mathbb{P}^n$}, i.e. a congruence, defined by the
parametric equations
\begin{align*}\begin{cases}
y_0&={\lambda},\\
y_i&={u^i\lambda-f^i({\mathbf u})\mu}, \ i=1,\dotsc,n-1\\
y_n&=\mu
\end{cases}\end{align*}
     where
\begin{itemize}
\item ${(\lambda:\mu)}\in\mathbb{P}^1$ are homogeneous coordinates on
a
     line of ${\mathcal B}$;
\item ${\mathbf u}=(u^1,\dotsc,u^{n-1})$
are the (local) parameters of ${\mathcal B}$;
\item $(y_0:\dotsb:y_n)$ are homogeneous  coordinates on
$\mathbb{P}^n$.
\end{itemize}

So for every ${\mathbf u}$ in the domain $\Omega$, the above equations
define a line $\Lambda({\mathbf u})$.

A dictionary can be written
translating properties of the system \eqref{eq:cons} to properties
of the family of lines, and conversely. For example, reciprocal
transformations of the system correspond to projectivities in
$\p^n$.
The eigenvalues $\lambda_1({\mathbf u}),\ldots ,
\lambda_{n-1}({\mathbf u})$ are in natural bijection with the foci of
${\mathcal B}$ on the line $\Lambda(\mathbf u)$. In particular, the system is
strictly hyperbolic if and only
if on a general line of ${\mathcal B}$ there are $n-1$ distinct foci
of the congruence.
The rarefaction curves correspond to developable ruled surfaces in $\p^n$.

\subsection{Temple Systems again}

In general for a hyperbolic system of conservation laws the corresponding
focal locus of the congruence ${\mathcal B}$ is a hypersurface $F$
and the lines of the family
are tangent to $F$ at the $n-1$ foci.
For a Temple system the situation is different.

\begin{thm} (\cite{AF})
The congruences of lines ${\mathcal B}$ associated to  Temple systems are
characterized by the properties that  every focus is a fundamental
point and that the
rarefaction curves correspond to planar pencils of lines of ${\mathcal B}$.
\end{thm}

By Theorem \ref{DC} it follows that dim $F=n-2$ and the order of
${\mathcal B}$ is one.

The classification of the $T$-systems up to reciprocal transformations
is equivalent to the classification of non necessarily algebraic
congruences of lines of
order one with planar pencils.

\begin{ex}\label{wave}
\emph{The wave equation}.
The equation $f_{tt}=f_{xx}$ can be rewritten as a system of two
conservation laws
\begin{align*} \begin{cases}
u^1_t=u^2_x \\
u^2_t=u^1_x
\end{cases}\end{align*}
The associated congruence ${\mathcal B}$ in $\p^3$ is
formed by the lines meeting two fixed
skew lines $L$ and $L'$. On each line $\Lambda$ of ${\mathcal B}$
there are two distinct foci, its
intersections with $L$ and $L'$.

This is also an example of a Temple system. The associated congruence
of lines is linear.
It is rather easy to prove (see Proposition~\ref{prop:lineari}) that 
all linear congruences, in any
projective space, are
associated to some
     Temple system.
Conversely:
\end{ex}

\begin{thm} (\cite{AF}, \cite{AF2})
All Temple systems in $2$ and $3$ variables give rise
to linear
congruences. In
particular they are all algebraic.
\end{thm}

Agafonov and Ferapontov have conjectured that  congruences of lines whose
developable surfaces are planar pencils of lines have always
algebraic focal varieties
(possibly reducible and singular).


\section{Linear Congruences of lines}\label{sec:4}
As we said in the introduction, a linear congruence of lines
${\mathcal B}$ in $\p^n$
has the form ${\mathcal B}=\Gr(1,n)\cap\Delta$, where $\Delta$ is a
linear subspace
of dimension $n\choose 2$ of $\p(\wedge^2 V)$, with
$V:=H^0(\OO_{\p^n}(1))^*$, the space
of the Pl\"ucker embedding of the Grassmannian.

The ${n+1}\choose 2$ coordinates of a line $\ell\in\Gr(1,n)$
can be interpreted as entries of a skew-symmetric
$(n+1)\times(n+1)$ matrix
$(p_{ij})_{i,j=0,\ldots,n}$ of rank $2$.  A hyperplane $H$ in
$\p(\wedge^2 V)$
has an equation of the form $\Sigma_{i,j=0}^n a_{ij}p_{ij}=0$,
and the $a_{ij}$'s are  coordinates of $H$ as a point in the
dual space $\p(\wedge^2 V^*)$.
Clearly, also the dual coordinates $a_{ij}$ can be interpreted
as entries of a skew-symmetric matrix $A$.
$\Delta$ is the intersection of the $n-1$ hyperplanes $H_1,\dotsc,H_{n-1}$
with equations:
\begin{align}\label{linear}
\sum_{i,j=0}^n a^1_{ij}p_{ij}&=0, &\dotsb&
   &\sum_{i,j=0}^n a^{n-1}_{ij}p_{ij}&=0
\end{align}
associated to matrices $A_1,\dotsc,A_{n-1}$.
In the dual space $H_1,\dotsc,H_{n-1}$  generate the dual
$(n-2)$-space $\check\Delta$.

Some general results about fundamental varieties of linear
congruences are given in \cite{BM}; in particular, it is proved that
the focal locus of a linear congruence ${\mathcal B}$
is the degeneracy locus $F$ of a morphism of
sheaves of the form
\begin{equation}\label{mf}
\phi:\OO_{\p^n}^{\oplus (n-1)}\rightarrow \Omega_{\p^n}(2).
\end{equation}

Explicitly, there is an isomorphism:
\begin{equation}\label{com}
H^0(\Omega_{\p^n}(2))\cong(\wedge^2 V)^*,
\end{equation}
   and so
a global section of $\Omega_{\p^n}(2)$ is
a skew-symmetric matrix
of type $(n+1)\times(n+1)$ with entries in the base field.
Then, the  morphism $\phi$  in \eqref{mf} is defined by the
$n-1$  skew-symmetric matrices $A_1,\dotsc,A_{n-1}$.
The corresponding
degeneracy locus $F$ in $\p^n$ is
defined by the equations
\begin{equation}\label{sks}
\sum_{i=1}^{n-1}\lambda_i A_i[X]=0
\end{equation}
for some $[\lambda]=(\lambda_1,\dotsc,\lambda_{n-1})\neq (0,\dotsc,0)$,
where $[X]$
denotes the column matrix of the coordinates.

Since the matrices $A_1,\dotsc,A_{n-1}$ are
skew-symmetric, the situation changes whether $n$ is even or odd.

\begin{prop}(\cite{D2})\label{prop:flc}
If $F$ is the focal locus of a  linear congruence ${\mathcal B}$ in
$\p^n$, then
\begin{enumerate}
\item  for ${\mathcal B}$ general,  $F$ is smooth if
$\dim(F)\le 3$;
\item if $n$ is even, for each $[\lambda]\in\p^{n-2}$ equation \eqref{sks}
has at least one solution, and
$F$ is  rational;
\item\label{odd} if $n$ is odd,
the vanishing of the Pfaffian of the matrix $\sum_{i=1}^{n-1}\lambda_i A_i$
   defines a
hypersurface $Z$ of degree $(n+1)/2$  in  $\p^{n-2}$
(in which $\lambda_1,\dotsc,\lambda_{n-1}$ are the coordinates).
Furthermore, if
$\phi$ is general, for a fixed point $[\lambda]\in Z$,
   equation \eqref{sks} has a line contained in $F$ as solution,
and $F$ results to be a scroll
over (an open set of) $Z$.
\end{enumerate}
Besides, in both cases
\begin{equation*}
\deg(F)=\frac{n^2-3n+4}{2}.
\end{equation*}
\end{prop}

Another known result, which can be deduced by the above description of
the focal locus of a linear congruence
and which explains the fact that linear congruences give Temple 
systems is the following:
\begin{prop}\label{prop:lineari}
Let $\mathcal B$ be 
a linear congruence in $\p^n$, with focal locus $F\subset\p^n$;
then $\mathcal B$ has order one and is the closure of the family of
the $(n-1)$-secant lines to $F$.
Moreover, if $P$ is a general point in $F$, the family of the lines 
of $\mathcal B$
through $P$ 
is a pencil, whose plane intersects $F$, out of $P$, in a curve of 
degree $n-2$.

In particular, a linear congruence corresponds to a Temple system.
\end{prop}

\begin{proof}
The first assertion follows from standard Schubert
calculus and the second one can be deduced from Theorems~\ref{DC} and 
\ref{prop:fofi}.
We give here a different direct proof of
both facts which relies on the previous observations.
 From Equation~\eqref{sks} (or
Equations~\eqref{linear})  we infer that through the general point 
$P\in \p^n$ with
coordinates $[X]=[x_0,\dotsc,x_n]$ there passes only the line of $\mathcal B$
whose coordinates are  the maximal 
minors of
\begin{equation}
\label{m1}
A:=
\begin{pmatrix}
\sum_i a^1_{0i}x_i&\ldots&\sum_i a^1_{n i}x_i  \\

\vdots&&\\
\sum_i a^{n-1}_{0i}x_i&&\sum_i a^{n-1}_{n i}x_i
\end{pmatrix}
.
     \end{equation}
Moreover, if we think of $A$ as a matrix with linear entries in $\p^n$, the
degeneracy locus of $A$ is $F$. Now, we can fix one line $\ell\in\mathcal B$
and without loss 
of generality
we can suppose that it does not intersect the $(n-2)$-dimensional space
defined by $x_0=x_1=0$; then $F\cap\ell$ is defined by the determinant of the
submatrix 
of $A$ formed by the last $(n-1)$ columns, and therefore---if $\ell$
is general---$F\cap\ell$ is a zero dimensional scheme of length $(n-1)$.

Now, if $P$ is a general focal point,
by the linearity the lines of $\mathcal B$ through it form a
pencil $\ell_P$: in fact in Equation~\eqref{sks} we can suppose that
$A_1[X]=0$ and $A_2[X]\neq,\dotsc, A_{n-1}[X]\neq 0$, \ie $P$ is the 
\emph{centre} 
  of
the linear complex $A_1$. Then the lines of 
$\ell_P$ are contained in the plane
whose Pl\"ucker coordinates are 
  given by the $(n-2)\times(n-2)$-minors of
\begin{equation*}
A':=
\begin{pmatrix}
\sum_i a^2_{0i}x_i&\ldots&\sum_i a^2_{n i}x_i  \\

\vdots&&\\
\sum_i a^{n-1}_{0i}x_i&&\sum_i a^{n-1}_{n i}x_i
\end{pmatrix}
.
     \end{equation*}
The plane $\pi_P$ of the 
pencil intersects $F$ in $P$ and in a curve
of degree $(n-2)$: in fact $\pi_P$ intersects the hypersurface $V$ defined by a
(fixed) minor of~\eqref{m1} in a curve of degree $(n-1)$  which splits in the
line of the congruence through $P$ and contained in 
$V$ and residually in a curve
which must be contained in $F$.
\end{proof}

The dual variety of the Grassmannian $\check{\Gr}(1,n)$ parametrises
the tangent
hyperplanes to $\Gr(1,n)$. It is defined by the maximal Pfaffians of
the matrix $A$,
therefore if $n$ is odd, it is the
hypersurface  in
$\check{\p}^{n}$ of degree $\frac{n+1}{2}$, defined by the Pfaffian,
while if $n$ is even, it has codimension $3$.
In the case $n$ odd, using coordinates
$(\lambda_1,\dotsc,\lambda_{n-1})$ in $\check\Delta$,  the intersection
$S:=\check{\Gr}(1,n)\cap\check{\Delta}$  is defined by
$\Pfaff(\sum_i\lambda_i A_i)=0$, hence it coincides
with the hypersurface $Z$ of Proposition \ref{prop:flc}.


\subsubsection{Linear Congruences in $\p^3$ and $\p^4$}
We end this section by briefly recalling the classification of the linear
congruences in low dimensional projective spaces.

In $\p^3$, the situation is very simple: $\Gr(1,3)\subset\p^5$ is the
Klein quadric,
its dual is again a
quadric, and $\check{\Delta}$ is a line. Correspondently, we have the
following cases: if $\check{\Delta}$ is general, it intersects
$\check{\Gr}(1,3)$ at
two distinct points and the congruence represents the join of the two
corresponding lines.
The line $\check{\Delta}$
can be tangent and therefore intersects $\check{\Gr}(1,3)$ in a double
point. Correspondently, we have a congruence which has as focal locus a double
line (and the congruence is a subset of the set of lines meeting its support).
Finally, it can happen that
$\check{\Delta}\subset\check{\Gr}(1,3)$:
in this case, we do not have a congruence, since the corresponding
family of lines in
$\p^3$ has dimension greater than three.

In $\p^4$, the situation is more complicated:  $\Gr(1,4)\subset\p^9$ has
dimension six and degree five. $\check{\Delta}\cong\p^2$ and in the
general  case
$\check{\Delta}\cap\check{\Gr}(1,4)=\emptyset$:  then the congruence
is given by the trisecants to a projected Veronese surface.
If $\check{\Delta}\cap\check{\Gr}(1,4)$ is non empty and finite, then
its length
is at most $3$. If it is a single point, this point gives us
a focal plane: the congruence is given by the secant lines to a cubic scroll
which meet this plane also. If $\check{\Delta}\cap\check{\Gr}(1,4)$ is
two points, then we have the lines meeting  two fixed
planes and which meet a quadric also. If
$\check{\Delta}\cap\check{\Gr}(1,4)$ is given by
three points,  we get the lines meeting  three
planes;  the focal locus has also a fourth component, which is parasitic,
the plane spanned by the $3$
points of intersection of  the three planes two by two. Of course
also  limit cases of these are possible, \eg if
$\check{\Delta}\cap\check{\Gr}(1,4)$ is a double point, etc.
Finally,  $\check{\Delta}\cap\check{\Gr}(1,4)$ can be a curve or
$\check{\Delta}\subset\check{\Gr}(1,4)$: several cases are possible
but they don't define  a congruence, since
the corresponding family of lines in
$\p^4$ has dimension greater than four (see the original paper of
Castelnuovo \cite{C}).

The article \cite{AF2} contains the interpretation of this 
classification in terms
of that of  Temple systems in three variables up to reciprocal transformations.

\section{Linear Congruences in $\p^5$}\label{sec:5}

A linear congruence of lines ${\mathcal B}$ in $\p^5$ is  of the form
${\mathcal B}=\mathbb G(1,5)\cap \Delta$,
where $\Delta$ is a linear space of dimension $10$,
intersection of $4$ hyperplanes
with equations:
\begin{align}
\sum_{i,j=0}^5 a_{ij}p_{ij}&=0, &\sum_{i,j=0}^5 b_{ij}p_{ij}&=0,&
\sum_{i,j=0}^5 c_{ij}p_{ij}&=0, &\sum_{i,j=0}^5 d_{ij}p_{ij}&=0
\end{align}
associated to skew-symmetric $6\times 6$  matrices $A,B,C,D$.

The dual variety of the Grassmannian  is the cubic hypersurface
   $\check{\Gr}(1,5)$ in
$\check{\p}^{14}$ defined by the Pfaffian.
One can think of
$\check{\Gr}(1,5)$ as the locus of skew-symmetric matrices (in the
$(a_{ij})$'s coordinates) of rank at most four.
Using coordinates
$(a,b,c,d)$ in $\check\Delta$,  the intersection
$S:=\check{\Gr}(1,5)\cap\check{\Delta}$  is defined by
$\Pfaff(aA+bB+cC+dD)=0$, hence, in general,
it is a cubic surface.

The rational Gauss map $\gamma: \check{\Gr}(1,5) \dasharrow \Gr(1,5)$
associates
to a tangent hyperplane its unique tangency point. It is regular
outside $\Sing(\check{\Gr}(1,5))$, which is
naturally isomorphic to $\Gr(3,5)$  and corresponds to the matrices of rank two
(which indeed is dually isomorphic to $\Gr(1,5)$ also).
If $A$ is such a matrix, the
corresponding $3$-space $\pi_A$ is the projectivised kernel of $A$.
As a hyperplane section of $\Gr(1,5)$,
$A$ represents the lines in $\p^5$ meeting $\pi_A$. In what follows
we will always identify $\Sing \check{\Gr}(1,5)$ with $\Gr(3,5)$.

For general $\Delta$, the image $\gamma(S)$ is a
$2$-dimensional family of lines, whose union is a smooth $3$-fold $F$
of degree $7$ and sectional genus $4$, called
Palatini scroll. The lines of the congruence ${\mathcal B}$ are the
$4$-secant lines of $F$ and $F$ is
the focal locus of ${\mathcal B}$ (see also \cite{O}, \cite{BM},
\cite{FM} and \cite{D2}). More in general:

\begin{prop}\label{prop:mio}
Let  $\Delta$ be a linear space of dimension $10$ in $\p^{14}$
     such that
     $S:=\check{\Delta}\cap\check{\Gr}(1,5)$ is a reduced  cubic surface
and $S \cap \
{\Gr}(3,5)$ is  empty. The lines $b\in\gamma(S)$ span
(set-theoretically) the focal divisor $F$ of the
linear congruence
${\mathcal B}=\mathbb G(1,5)\cap \Delta$, \ie  the map $f$
of \textup{\eqref{eq:f}} drops rank exactly at the pairs $(b,P)$ such
that $b\in\gamma(S)$ and $P\in\Lambda(b)$.
The lines of ${\mathcal B}$ are the $4$-secants of $F$.
\end{prop}

\begin{proof}
Since $\gamma(S)\subset {\mathcal B}$ is the
set of points whose tangent hyperplanes (to $\Gr(1,5)$) contain
$\Delta(=\gen{{\mathcal B}})$, we have one
inclusion recalling that  for these points the global characteristic
map---see for example Definition~1 of
\cite{DP2}---of the congruence drops rank.  For the other inclusion,
given a point $P\in F\cap\Lambda(b)$, by
definition, through a smooth
$b\in {\mathcal B}$ there passes at least one tangent
direction to $\Gr(1,5)\setminus {\mathcal B}$ contained in $\Delta$:
we can find---by dimensional reasons---a tangent hyperplane of $P$
containing $\Delta$ also.  The last assertion
follows from Theorem \ref{prop:fofi}.
     \end{proof}

Classifying linear congruences in $\p^5$ amounts to describing all
special positions of the
$3$-space $\check\Delta$ with respect to $\check{\Gr}(1,5)$ and to its singular
locus.
As we will see the situation is rather complicated and several
different cases are possible.

\begin{rmk}
Actually, we are interested mainly to the case of ``true'' linear congruences,
\ie
if the intersection ${\mathcal B}=\Delta\cap\Gr(1,5)$
is proper  and the focal locus has dimension three.
Therefore, in what follows we will not give many details on the
cases such that
${\mathcal B}$  has dimension $>4$ or
$F$ has dimension $>3$.
\end{rmk}

\begin{prop}
With the above notation, if $\check{\Delta}$ is contained in
$\check{\Gr}(1,5)$, then $\check{\Delta}$ meets $\Gr(3,5)$ and
     the focal locus $F$ of ${\mathcal B}$ has dimension $>3$.
If $\check{\Delta}\cap\check{\Gr}(1,5)$ is a cubic surface $S$, then
$\dim F>3$ if and only if
$S$  intersects $\Gr(3,5)$  at least along a curve.
\end{prop}

\begin{proof}
The first claim follows from Corollary 11 of \cite{MM}. The second
and the third one are
     Theorem 4.3 and Theorem 4.9 of
\cite{FM}.
\end{proof}

    From now on, we will always assume that the congruence ${\mathcal
B}$ is obtained from a
10-space $\Delta$ such that the intersection
$\check{\Delta}\cap\check{\Gr}(1,5)$ is proper \ie a cubic surface $S$.
    From \cite{FM}, Remark~4.4, it follows that every cubic surface in
$\p^3$ can be realized in this way.

If $S$ is smooth, then $F$ is always a Palatini scroll. For any $S$,
possibly singular
and/or reducible, such that $\dim F=3$, the Hilbert polynomial of $F$
is $P_F(t)=7/6 t^3+ 2t^2 + 11/6t + 1$. The equations of $F$ can be written
explicitly as  maximal minors of the following $4\times 6$ matrix
(see \cite{FM}):
\begin{equation}
\label{matrix}
M=
\begin{pmatrix}
\sum_i a_{0i}x_i&.&.&\sum_i a_{5i}x_i  \\
\sum_i b_{0i}x_i&.&.&\sum_i b_{5i}x_i \\
\sum_i c_{0i}x_i&.&.&\sum_i c_{5i}x_i\\
\sum_i d_{0i}x_i&.&.&\sum_i d_{5i}x_i
\end{pmatrix}
     \end{equation}
that we will also write in the form
\begin{equation}
\label{matrix2}
\begin{pmatrix}
L_{10}&L_{11}&.&.&.&L_{15}  \\
L_{20}&L_{21}&.&.&.&L_{25} \\
L_{30}&L_{31}&.&.&.&L_{35}\\
L_{40}&L_{41}&.&.&.&L_{45}
\end{pmatrix}
     \end{equation}

The lines of ${\mathcal B}$ through a general point $P$ in $F$ form a
linear pencil
contained in a plane cutting $F$, out of $P$, along a plane cubic curve.
The coefficients of the equations of this plane are  the lines of the
matrix \eqref{matrix}
computed at $P$.

\subsubsection*{}
It is well known
that a cubic surface with isolated singularities can have at most $4$
double points,
or one triple point if it is a cone, and that an irreducible cubic
with a singular curve is necessarily ruled with a double line.
Observe also that if
$A$ is a singular point of
$S$, then either
$A\in\Gr(3,5)$ or
$\check\Delta$ is tangent to
$\check{\Gr}(1,5)$ at $A$. We will now
consider the various possibilities for $S$, studying the linear congruences
${\mathcal B}$ and their focal loci. We will also give some explicit
examples, constructed using \cocoa
(see  \cite{cocoa}).

We begin by studying two special situations for the singularities
of $S$, \ie the cases when one or more singular points belong to
$\Gr(3,5)$ and when $S$ is reducible  not meeting $\Gr(3,5)$.

\subsection{$S$ with singular points on $\Gr(3,5)$}
\subsubsection{Only one singular point} Let $A$ be a singular point
of $S$
belonging to
$\Gr(3,5)$. In this case $\Delta$
is contained in the linear complex of the lines meeting the $3$-space
$\pi_A$,
hence each line of
the congruence ${\mathcal B}$ intersects $\pi_A$. We get that $F$ splits as
$F=\pi_A\cup Y$, where $Y$ is an arithmetically Cohen-Macaulay (aCM
for short) threefold
of degree $6$ and sectional genus $3$, defined by the maximal minors
of a $3\times 4$ matrix
of linear forms.
Indeed, if we choose a system of coordinates such that $\pi_A$ has
equations $x_0=x_1=0$, then
the only non-zero coordinate of $A$ is
$a_{01}$ and  the matrix \eqref{matrix2}
becomes
\begin{equation}
\label{matrix3}
\begin{pmatrix}
x_1&-x_{0}&0&0&0&0  \\
L_{20}&L_{21}&.&.&.&L_{25} \\
L_{30}&L_{31}&.&.&.&L_{35}\\
L_{40}&L_{41}&.&.&.&L_{45}
\end{pmatrix}.
     \end{equation}

One checks that $\pi_A$ is a component of $F$
and the residual is defined by the $3\times 3$ minors of the matrix
formed by the last
$3$ rows and $4$ columns of \eqref{matrix3}.

If the other generators of $\check\Delta$ are general, then $S$ is
smooth outside $A$. $Y$ results to be a
singular Bordiga scroll, with $6$ singular points all belonging to
$Y\cap \pi_A$, which  is a smooth
quadric surface.
The foci on a general line of  the congruence are contained one in
$\pi_A$ and the other three in $Y$, so
the lines of ${\mathcal B}$ are the trisecants of $Y$
meeting also $\pi_A$. The lines of ${\mathcal B}$ through a point of
$\pi_A$ cut  $Y$ at the points of a  plane
cubic, whereas  those through a point $P$ in $Y$ cut $Y$ residually
along a conic $C$ and $\pi_A$ along a line $L$.

\subsubsection{Two singular points}
Assume now that $S$ has also a second singular point $B$ on
$\Gr(3,5)$. Then $F$ is of the form
$F=\pi_A\cup \pi_B\cup Z$, where $Z$ is an aCM threefold of degree 5
and sectional genus 2.
To see this, let us note that the 3-spaces $\pi_A$ and $\pi_B$ are
     in general position (otherwise the line $\overline{AB}$ would be
contained in $\Gr(3,5)$), so we can assume that
$\pi_A$ is defined by $x_0=x_1=0$ and $\pi_B$ by $x_2=x_3=0$. So the
matrix $M$ becomes
\begin{equation}
\label{matrix4}
\begin{pmatrix}
x_1&-x_{0}&0&0&0&0  \\
0&0&x_3&-x_2&0&0 \\
L_{30}&L_{31}&.&.&.&L_{35}\\
L_{40}&L_{41}&.&.&.&L_{45}
\end{pmatrix}.
     \end{equation}
By developing its $4\times 4$ minors, we get that the two 3-spaces
$\pi_A$ and $\pi_B$ are irreducible components
of $F$ and the residual is defined by the $2\times 2$ minors of
\begin{equation*}
\begin{pmatrix}
x_0L_{30}+x_1L_{31}&L_{34}&L_{35} \\
x_0L_{40}+x_1L_{41}&L_{44}&L_{45}
\end{pmatrix},
     \end{equation*}
or, equivalently, of
\begin{equation*}
\begin{pmatrix}
x_2L_{32}+x_3L_{33}&L_{34}&L_{35} \\
x_2L_{42}+x_3L_{43}&L_{44}&L_{45}
\end{pmatrix}.
     \end{equation*}
$Z$ is therefore a  Castelnuovo threefold.  If the other generators
of $\check\Delta$ are general,
then $S$ is smooth outside $A$ and $B$. $Z$ results to be
singular at 8 points, 4 of them belonging to $Q_A:=Z\cap \pi_A$, and
the other 4 to
$Q_B:=Z\cap \pi_B$. $Q_A$ and $Q_B$ are both quadric surfaces, linear
sections of the unique
quadric containing $Z$, of equation $L_{34}L_{45}-L_{35}L_{44}=0$.
They intersect along the line
$\pi_A\cap\pi_B$.
Clearly the lines of ${\mathcal B}$ are the secants of $Z$
meeting also $\pi_A$ and $\pi_B$.
The lines of ${\mathcal B}$ through a point in $\pi_A$ or $\pi_B$
cut  $Z$ along the points of a conic, whereas
those through a point $P$ in $Z$ cut $Z$, $\pi_A$ and $\pi_B$
residually each along a line.

\subsubsection{Three singular points}
We assume now that $S$ has also a third singular point $C$
on $\Gr(3,5)$. Then $F$ is of the form
$F=\pi_A\cup \pi_B\cup \pi_C\cup V$, where $V$ is a
complete intersection of two quadrics.
Indeed, as in the previous case
     the 3-spaces $\pi_A$, $\pi_B$ and $\pi_C$ are two by two
     in general position, so we can assume that
$\pi_A$ is defined by $x_0=x_1=0$, $\pi_B$ by $x_2=x_3=0$
and $\pi_C$ by $x_4=x_5=0$.
The matrix $M$ takes the form
\begin{equation*}
\begin{pmatrix}
x_1&-x_{0}&0&0&0&0  \\
0&0&x_3&-x_2&0&0 \\
0&0&0&0&x_5&-x_4\\
L_{0}&.&.&.&.&L_{5}
\end{pmatrix}.
     \end{equation*}
It is easy to see that the residual of the 3 spaces
$\pi_A$, $\pi_B$ and $\pi_C$
is defined by the equations $x_0L_0+x_1L_1=x_4L_4+x_5L_5=0$,
or equivalently by
$x_2L_2+x_3L_3=x_5L_5+x_4L_4=0$.
$V$ is therefore a  Del Pezzo threefold.  Again, if the
other generators of $\check\Delta$ are general,
then $S$ is smooth outside $A$, $B$ and $C$. Also $V$ results to be
smooth. The intersections of $V$ with the spaces $\pi_A$,
$\pi_B$ and $\pi_C$ are quadric surfaces
intersecting two by two  along
lines.
The lines of ${\mathcal B}$ are the lines meeting
simultaneously $\pi_A$, $\pi_B$, $\pi_C$ and $V$.

\subsubsection{Four singular points}
We consider finally the case in which $S$ has also a
fourth singular point $D$ on $\Gr(3,5)$. Then $F$ is of the
form
$F=\pi_A\cup \pi_B\cup \pi_C\cup \pi_D \cup W$, where
$W$ is a rational normal cubic scroll.

It is clear that  the four spaces $\pi_A$, $\pi_B$, $\pi_C$
and $\pi_D$ are all components of the focal locus $F$,
and that the lines of ${\mathcal B}$ are characterized
by the property of meeting simultaneously all of them.
The component $W$, that must exist by degree reasons,
arises in the following way. Let $H$ be a hyperplane
containing $\pi_A$, it intersects the other three spaces
along  planes, call them $\alpha_B$, $\alpha_C$ and
$\alpha_D$,  and ${\mathcal B}_H$, the restriction of
${\mathcal B}$ to $H$, is formed by the lines meeting them.
There is a uniquely determined fourth  plane meeting
$\alpha_B$, $\alpha_C$ and
$\alpha_D$ along lines, and it is necessarily contained
in the focal locus of ${\mathcal B}_H$ (it is a
\lq\lq parasitic plane'', see \cite{i}, because  its lines intersect
all the three planes $\alpha_B$, $\alpha_C$
and
$\alpha_D$ but a general line of ${\mathcal B}_H$ does not meet it).
If we let $H$ vary, we get a 1-dimensional
family of planes, whose union is $W$.

\subsubsection{Special degenerate cases}
These 4 examples don't exhaust the list of possibilities for the
surfaces $S$ meeting $\Gr(3,5)$. For instance the space
$\check\Delta$ could be contained in the tangent space to
$\Gr(3,5)$ at one of the points of intersection, or  intersect it
along a plane or a line. In this situation four, three or two of
the points of intersection of $\check\Delta$ with $\Gr(3,5)$ get
identified. Hence some component of the focal locus appears with
multiplicity greater than one and the 4 foci on a general line of ${\mathcal
B}$ are not distinct. A particular case is when $S$ is a cone of
vertex a point $A$ of $\Gr(3,5)$, then the 3-space corresponding to
the vertex counts twice as component of $F$.

\subsection{$S$ with a double line not meeting $\Gr(3,5)$} \label{rid}
Assume that $S$ is reducible and disjoint from
$\Gr(3,5)$, therefore of the form $S=\pi\cup Q$ where $\pi$ is a
plane and $Q$ a (possibly reducible) quadric.

In \cite{MM} the linear spaces contained in $\check\Gr(1,5)$ and not
meeting $\Gr(3,5)$ are classified up to the natural action of
$\PGL_6$. They can be interpreted also as linear spaces of
skew-symmetric matrices of constant rank 4.

It results that
there are two orbits of lines, an open irreducible orbit of dimension 22,
and a codimension one closed orbit. Note that the
Gauss map $\gamma$, if restricted to such a line $\ell$, is
regular and defined by the derivatives of the cubic Pfaffian
polynomial (the equation of $\check\Gr(1,5)$), which have degree
$2$. Therefore $\gamma(\ell)$ is embedded in $\p^{14}$ as a conic.
The lines parametrised by this conic represent a ruling of a smooth quadric
or the the lines of a quadric cone, respectively for the two orbits.

The planes are distributed in 4 orbits, all of dimension 26.
The Gauss image of such a plane $\pi$  is embedded in
$\p^{14}$ as a Veronese surface $V$. The four orbits correspond precisely
to the possible
double Veronese embeddings of the plane in $\Gr(1,5)$.
The lines represented by $V$ have the following
geometrical interpretation, in the four cases (see \cite{SU}):
\begin{enumerate}
      \item\label{1} the secant lines of a skew cubic curve $\mathcal C$
      embedded in a 3-space $L$;
      \item\label{2} the lines contained in a quadric 3-fold $\Gamma$
and meeting a
      fixed line $r$ contained in $\Gamma$;
      \item\label{3} the lines joining the corresponding points in a fixed
      isomorphism between two disjoint planes;
      \item\label{4} the lines of a cone $C(\mathcal V)$ over a
projected Veronese
      surface.
\end{enumerate}

If $\check\Delta$ is generated by one of these four types of
planes plus one general point of $\check\p^{14}$, then the
residual component of $\pi$ in $S$ is a quadric surface of maximal
rank. We describe now the corresponding congruences and their
focal loci.

\subsubsection{Case~\eqref{1}}  $F= L'\cup X$, where $L'$ is a scheme
whose support is $L$,
      having $\mathcal C$ as embedded component, and $X$ is a singular
threefold of degree 6
      with Hilbert polynomial $P_X(t)=t^3+3t^2+2$; $X$ meets $L$ along a
      quartic surface with $\mathcal C$ as singular locus. Through a general
      point of $L$ there passes a 1-dimensional family of lines of
      $\mathcal B$,  all contained in $L$, and through a general point of
$\mathcal C$
      a 2-dimensional family.  $L$ and $\mathcal C$ are components of the
focal locus,
      but a general line of $\mathcal B$ does not meet $L$ nor $\mathcal C$,
      they are therefore both
      parasitic components. On a general line of $\mathcal B$ the
       4 foci lie all on $X$. It follows that $X$ is an example of
       (singular) threefold containing a
       3-dimensional family of plane cubic curves (see \cite{MP}).

       An explicit example of this type of congruence is given by the following
       skew-symmetric matrix, which represents $\check\Delta$:
       \begin{equation*}
     \begin{pmatrix}
         0 & -d & a & b & c & 0 \\
         . & 0 & 0 & a & b & c \\
         . & . & 0 & d & 0 & d \\
         . & . & . & 0 & 0 & 0 \\
         . & . & . & . & 0 & -d \\
         . & . & . & . & . & 0 \\
       \end{pmatrix}.
\end{equation*}

\subsubsection{Case~\eqref{2}}  $F=\Gamma '\cup Y$, where $\Gamma'$ is a scheme
      whose support is $\Gamma$,
      having $r$ as embedded component, and $Y$ is a quintic
      threefold with Hilbert polynomial $P_Y(t)=5/6t^3+5/2t^2+5/3t+1$.
      Through a general point in $\Gamma$ there passes a $1$-dimensional
      family of lines, through a general point in $r$ a $2$-dimensional
      family. Hence $r$ is a parasitic scheme. $Y$ meets $\Gamma$ in a
quartic surface
which is singular along $r$.

An explicit example is provided by the matrix:
     \begin{equation*}
     \begin{pmatrix}
         0 & a & b & c & d & 0 \\
         . & 0 & 0 & d & c & b \\
         . & . & 0 & 0 & d & 0 \\
         . & . & . & 0 & 0 & d \\
         . & . & . & . & 0 & a \\
         . & . & . & . & . & 0 \\
       \end{pmatrix}.
\end{equation*}

\subsubsection{Case~\eqref{3}}  $F=Z_1\cup Z_2$, where $Z_1$ is a
rational cubic
      scroll, \ie $\p^1\times \p^2$, which is the union of the lines
      parametrised by $\check\Delta$,
       and $Z_2$ is a singular Del Pezzo
      threefold, complete intersection of two quadrics. $Z_1\cap Z_2$
      is a quartic surface, rational normal scroll of type $(2,2)$,
      which cuts a conic on each plane of $Z_1$.
      The singular locus of $Z_2$ is the union of two conics.
      The foci of a general
      line lie two on each component of $F$, so $\mathcal B$ is formed
      by the lines  which are secant both $Z_1$ and $Z_2$.

     An  example of this congruence is the following:
        \begin{equation*}
     \begin{pmatrix}
         0 & a & b & d & 0 & 0 \\
         . & 0 & c & 0 & d & 0 \\
         . & . & 0 & 0 & 0 & d \\
         . & . & . & 0 & a & b \\
         . & . & . & . & 0 & c \\
         . & . & . & . & . & 0 \\
       \end{pmatrix}.
\end{equation*}

\subsubsection{Case~\eqref{4}} $F=C(\mathcal V)\cup T$, where $T$ is
a cubic threefold
contained in a hyperplane of $\p^5$ passing through the vertex of the cone.
The lines of $\mathcal B$ are the trisecants of $F=C(\mathcal V)$,
meeting also $T$.

     An  example of this congruence is the following:
        \begin{equation*}
    \begin{pmatrix}
         0 & a & b & c & d & 0 \\
         . & 0 & c & d & b & 0 \\
         . & . & 0 & 0 & 0 & d \\
         . & . & . & 0 & a & 0 \\
         . & . & . & . & 0 & 0 \\
         . & . & . & . & . & 0 \\
       \end{pmatrix}.
\end{equation*}

\subsubsection{Special degenerate cases}
In all the four cases, for special choices of the $4^\textup{th}$
generator of $\check\Delta$, the rank of the quadric $Q$ can
decrease, giving raise to various degenerations of the focal loci
of the just considered congruences.

\subsection{Classification of linear congruences in $\p^5$.}
Table~\ref{t:1} summarizes the sketch of classification of
linear congruences of lines in $\p^5$ with $3$-dimensional focal
locus.

Note that the smooth Palatini threefolds, which are scroll over cubic
surfaces, have a non-trivial moduli space. Indeed the moduli space of cubic
surfaces is rational of dimension $4$, and the moduli space of rank
$2$ vector
bundles $E$ on a cubic surface $S$, such that $\mathbb P(E)$ is a Palatini
scroll, has dimension $5$ (see \cite{FM}).

We finish this article with a few words about the singular Palatini
scrolls and the corresponding congruences.
Actually, the singular Palatini scrolls which come out from  a surface $S$
with (isolated) double
points $P_i$ (such that in fact $\Delta$ is tangent to
$\check(\Gr(1,5)\setminus\Gr(3,5)$) at $P_i$)
have only isolated singularities and can be described as the smooth
ones. Therefore also the congruences associated to these
can be described as in the case of the smooth Palatini scroll.

Let us pass now to the two more problematic Palatini scrolls, \ie the
ones which come from an $S$
that is either a cubic cone or a ruled surface, singular along a
line. In both cases, $S$ has a $1$-dimensional
family of lines, and these lines correspond to quadric surfaces in $\p^5$.

A $3$-space $\check\Delta$ such that $\check\Delta\cap\check\Gr(1,5)$
is a cubic cone of
vertex $A$ must be contained in the quadric $\ii_A$, the second
fundamental form of $\check\Gr(1,5)$ at
$A$, and  intersect the vertex of $\ii_A$ only at $A$. An example is
the following:
       \begin{equation*}
     \begin{pmatrix}
         0 & c & a & d-c & b & 0 \\
         . & 0 & d & a & 0 & b \\
         . & . & 0 & c & d & 0 \\
         . & . & . & 0 & 0 & c \\
         . & . & . & . & 0 & 0 \\
         . & . & . & . & . & 0 \\
       \end{pmatrix}.
\end{equation*}
The Palatini threefolds obtained in this way result to be singular
along the line corresponding to $A$. This Palatini scroll can indeed
be constructed
as the smooth one, and also the congruence does not give problems.

To obtain a ruled cubic $S$ with double line $\ell$, disjoint from
$\Gr(3,5)$, $\check\Delta$ must be
contained in the intersection of the tangent spaces to
$\check\Gr(1,5)$ at $P$, for all $P\in r$. One
checks that such a $3$-space exists only if $\ell$ is a line of the
closed orbit (see Subsection~\ref{rid}),
whereas  in the other case
no $3$-space is contained in the required intersection.
The Palatini scroll constructed in such a way is singular along the
quadric cone $Q$ which correspond
to $\ell$. Now, if we consider the secant lines to $F$ which meet
also $Q$ form---by dimensional reasons---a congruence.
Therefore this is the congruence which we are looking for,
since on a general line of
it we have two distinct foci plus a double focus.

\begin{table}
\begin{center}
\textup{
\begin{tabular}{|m{3.5cm}|m{2.2cm}|m{3cm}|m{3cm}|}
\hline
Surface $S$
    & $\Sing S$
     & Focal locus $F$
      & Remarks
       \\
\hline
\hline
smooth
    &
     & smooth Palatini scroll
      & $F$ has moduli
       \\
\hline
irreducible with isolated singularities $P_1,\dotsc , P_k$, $k\leq 4$
    & $P_i\notin\Gr(3,5)$
     & singular Palatini scroll
      & $F$ can be singular at a point or along a line
       \\
\cline{2-4}
    & $P_i\in\Gr(3,5)$
     & the $3$-space $\pi_{P_i}$ is an irreducible component of $F$
      & the residual intersects $\pi_{P_i}$ along a smooth quadric
       \\
\hline
irreducible ruled with a double line $r$
    & $r\cap\Gr(3,5)=\emptyset$
     & singular Palatini scroll
      & $F$ is singular along a quadric cone
        $\check\Delta$  is tangent to $\check\Gr(1,5)$ along $r$
       \\
\hline
reducible $\pi\cup Q$, $\pi$ a plane, $Q$ a quadric
    & conic in $\pi$,  disjoint from $\Gr(3,5)$
     & $L'\cup X$, $\deg X=6$
      & $L'$ is a parasitic component
       \\
\cline{3-4}
    &
     & $\Gamma'\cup Y$, $\deg Y=5$
      & $\Gamma'$ has a parasitic line as embedded component
       \\
\cline{3-4}
    &
     & $Z_1\cup Z_2$, $\deg Z_2=4$, $Z_1~=~\p^1~\times~\p^2$
      & the lines of $\mathcal B$ are secant both $Z_1$ and $Z_2$
       \\
\cline{3-4}
    &
     & $C(\mathcal V)\cup T$, $\deg T=3$
      & the lines of $\mathcal B$ are the trisecants $C(\mathcal V)$ meeting $T$
       \\
\hline
successive degenerations of the previous cases
    &
     &
      &
       \\
\hline
\end{tabular}}

\end{center}
\caption{}\label{t:1}
\end{table}

\subsubsection*{}
Note that this is only a sketch of classification. The task of
writing a complete classification of
linear congruences up to  isomorphism seems to us out of reach,
because of the multitude
of special cases, due to the moduli of cubic surfaces (see for
instance \cite {BL}) and to their several
different embeddings in
$\check\Gr(1,5)$.

In a forthcoming paper we plan to apply these results to the
classification of Temple systems
in 4 variables. We have  to point out that this classification holds
on an algebraically
closed field, so it will be necessary to refine it over the real field.

     \providecommand{\bysame}{\leavevmode\hbox
to3em{\hrulefill}\thinspace}

\end{document}